\documentclass[aap,MSNbibl,citesort,dvips]{arximspdf}
\usepackage{mathbh}
\usepackage{accents}

\doi{10.1214/10-AAP703}
\volume{21}
\issue{4}
\pubyear{2011}
\firstpage{1253}
\lastpage{1281}

\makeatletter

\newtheorem{theorem}{Theorem}[section]
\newtheorem{corollary}[theorem]{Corollary}
\newtheorem{lemma}[theorem]{Lemma}
\newtheorem{proposition}[theorem]{Proposition}

\newproclaim{definition}[theorem]{Definition}
\newproclaim{example}[theorem]{Example}
\newproclaim{remark}[theorem]{Remark}

\newcommand{\vect}[1]{\mathbf{#1}}

\newcommand{\Ii}{\mathbh{1}}
\newcommand{\iLn}{\accentset{\circ}{\Lambda}_n}
\newcommand{\iLnp}{\accentset{\circ}{\Lambda}_{n+1}}

\newcommand{\NN}{\mathbb{N}}
\newcommand{\RR}{\mathbb{R}}
\newcommand{\ZZ}{\mathbb{Z}}
\newcommand{\EE}{\mathbb{E}}
\newcommand{\N}{{\mathcal N}}
\newcommand{\Ns}{{\mathcal N_{\mathrm{supp}}}}
\newcommand{\C}{{\mathcal C}}
\newcommand{\F}{{\mathcal F}}
\newcommand{\M}{{\mathcal M}}
\newcommand{\x}{\mathbf{x}}
\newcommand{\y}{\mathbf{y}}

\newcommand{\LN}{{\Lambda_N}}
\newcommand{\E}[1]{\langle #1\rangle}
\newcommand{\Var}{\operatorname{Var}}
\newcommand{\supp}{\operatorname{supp}}

\makeatother

\begin{document}
\begin{frontmatter}

\title{Necessary and sufficient conditions for realizability of point
processes}
\runtitle{Conditions for realizability}

\begin{aug}
\author[A]{\fnms{Tobias} \snm{Kuna}\corref{}\thanksref{t1,t2}\ead[label=e1]{t.kuna@reading.ac.uk}},
\author[B]{\fnms{Joel L.} \snm{Lebowitz}\thanksref{t1}\ead[label=e2]{lebowitz@math.rutgers.edu}}
and
\author[B]{\fnms{Eugene R.} \snm{Speer}\ead[label=e3]{speer@math.rutgers.edu}}
\runauthor{T. Kuna, J. L. Lebowitz and E. R. Speer}
\affiliation{University of Reading, Rutgers University and Rutgers University}
\address[A]{T. Kuna\\
Department of Mathematics\\
University of Reading\\
Whiteknights, P.O. Box 220\\
Reading RG6 6AX\\
United Kingdom\\
\printead{e1}}
\address[B]{J. L. Lebowitz\\
E. R. Speer\\
Department of Mathematics\\
Rutgers University\\
New Brunswick, New Jersey 08903\\
USA\\
\printead{e2}\\
\phantom{E-mail: }\printead*{e3}}
\end{aug}

\thankstext{t1}{Supported by NSF Grant DMR-0442066 and AFOSR Grant
AF-FA9550-04.}
\thankstext{t2}{Supported by Feodor Lynen fellowship of the A. von
Humboldt Foundation.}

\received{\smonth{11} \syear{2009}}

%
\begin{abstract}
We give necessary and sufficient conditions for a pair of
(generali\-zed) functions $\rho_1(\mathbf{r}_1)$ and $\rho_2(\mathbf{r}_1,\mathbf{r}_2)$,
$\mathbf{r}_i\in X$, to be the density and pair correlations of some point process
in a topological space~$X$, for example, $\RR^d$, $\ZZ^d$ or a subset
of these.
This is an infinite-dimensional version of the classical ``truncated
moment'' problem. Standard techniques apply in the case in which there can
be only a~bounded number of points in any compact subset of $X$. Without
this restriction we obtain, for compact $X$, strengthened conditions which
are neces\-sary and sufficient for the existence of a process satisfying a
further~re\-quirement---the existence of a finite third order moment. We
generalize the latter conditions in two distinct ways when $X$ is not
compact.
\vspace*{3pt}
\end{abstract}

%

\end{frontmatter}

\section{Introduction}\label{intro}

A \textit{point process} is a probability measure on the family of all locally
finite configurations of points in some topological space $X$; for an
overview see \cite{DV-J}. Here we will often adopt the terminology of
statistical mechanics, referring to the points as \textit{particles} and to
their expected densities and correlations as \textit{correlation functions}.
In many applications, quantities of interest can be calculated from the
first few correlation functions---often the first two---alone (see
\cite{Hansen} and below). Given the process in some explicit form, for
example, as a Gibbs measure, one can in principle calculate these
correlation functions, although in practice this is often impossible. On
the other hand, one may start with certain prescribed correlation
functions; these might arise as computable approximations to those of some
computationally intractable process as occurs in the study of equilibrium
fluids \cite{Hansen} or might express some partial information about
an as
yet unknown process as in the study of heterogeneous materials. One would
like to determine whether or not these given functions are in fact the
correlation functions of \textit{some} point process, that is, are they
\textit{realizable}?

This paper is a continuation of our previous work on the realizability
problem \cite{KLS,CKLS} to which the reader may wish to refer but is
independent and may be read separately. In this Introduction we briefly
summarize our approach to the problem and then discuss a few
applications in the physical and biological sciences. We summarize
definitions and background in Sec\-tion~\ref{ptproc} and describe our new
results formally in the remainder of the paper.

It is often convenient to view the realizability problem as a
\textit{truncated moment problem}. In that setting it is
an instance of the general problem
of finding a process supported on some given subset of a linear space and
having specified \textit{moments} up to some given order, that is,
specified expectation values of all linear functions and their products up
to that order or equivalently of all polynomial functions of that degree.
(The problem is called ``truncated'' because not all moments are
prescribed.) To identify the realizability problem as a truncated moment
problem we use the interpretation of a configuration of points as a sum
of Dirac point measures and thus as a Radon measure on $X$. In
this sense the set of all point configurations becomes a subset of the
linear space of all signed Radon measures and to specify the correlation
functions of the process up to some order $n$ is then just to specify
moments, in the above sense, up to order $n$. This is an
infinite-dimensional instance of a classical mathematical problem \cite{AhKr62,trunc}.
For the one-dimensional moment problem there are many
powerful and interesting results but for higher-dimensional truncated
problems there are fewer (see \cite{trunc,Multi,CuFi08} and references
therein). For the infinite-dimensional truncated moment problem
we are only aware of \cite{Us73}.

In this paper we derive several classes of conditions on correlation
functions which are necessary and/or sufficient for their realizability by
a point process (or, in some cases, by a point process with certain extra
properties); for simplicity we suppose that we wish to realize only two
moments, that is, the first and second correlation functions but the methods
extend directly to the general case. The conditions we consider are
obtained via a standard general technique for moment problems, that is,
Riesz' method
\cite{Ak65}: one defines a linear functional on the space of quadratic
polynomial functions of configurations in such a way that its value on any
polynomial coincides with the expectation of the polynomial with
respect to
any realizing measure (should one exist). If the polynomial in question is
nonnegative on the set of point configurations then it is a necessary
condition for realizability that the value of the linear functional on the
polynomial be nonnegative as well. The linear functional is
expressible in
terms of the prescribed correlation functions alone so that this gives
rise to a necessary condition for realizability. These necessary conditions
are discussed in Section \ref{subsecnec}.

The challenge is to show that these conditions are in fact also sufficient
or, if they are not, to find an appropriate strengthening. There are two
general classical approaches to the construction of realizing
measures, one based on dual cones (see \cite{trunc}, Chapters I and V)
and one on an extension theorem for
nonnegative functionals (see Riesz's method in \cite{Ak65}).
We follow the latter path: we first extend the
linear functional described above to an appropriate larger space of
continuous functions then prove that the extended functional can be
realized by a measure. The first step follows in great generality from the
Riesz--Krein extension theorem, cf. Theorem \ref{RKET}.

If the set of
particle configurations is compact then the second step can be established
by the well-known Riesz--Markov theorem; in Section \ref{susechard} we use
these ideas to establish sufficiency of the conditions described above in
this case. The set of configurations is compact if the system has a local
restriction on the number of particles; such a restriction can arise
naturally from an a priori restriction on the total number of particles in
the system (the result in this case was already proven directly in
\cite{Percus2} and \cite{GarrodPercus}) or in a~setting where the given
correlation functions, by vanishing on certain sets (as would be
implied by
a hard-core exclusion condition), prohibit particles from being closer to
each other than some given distance. Lattice systems in which there can be
no more than a given number of particles per site are included in this case
(see Section \ref{susechard}).

When it is not known that the support of the desired measure is compact,
we use a compact function, that is, a function with compact level sets
to obtain an analogue of the Riesz--Markov representation theorem from the
Daniell theory of integration. In our case we may use an appropriate power
of a linear function as the desired compact function. In general, however,
the process obtained in this way will not automatically realize the highest
prescribed moment (cf. \cite{trunc}, Chapter V.1); this is a feature of
the truncated moment problem on noncompact spaces (in our case we must
even consider nonlocally compact spaces) which does not arise if
moments of all orders are prescribed because there is always a higher
moment at hand to control the lower ones. This difficulty is not avoidable;
in fact the conditions derived from the positive polynomials are not in
general sufficient in the noncompact case (see \cite{trunc} and
Example \ref{funny} below). An alternative approach for the locally compact
case is given in \cite{trunc}.

In this paper we propose a new and quite natural approach for
infinite-dimensional moment problems. We modify the
conditions in order that they become sufficient but they then
cease to be necessary; rather, they are necessary and sufficient for the
modified realizability problem in which one requires realizability of the
first two correlation functions by a process which has a finite third
(local) moment. This, to our knowledge, is the first extension of the abstract
characterization of necessary and in some sense sufficient conditions
for an
infinite-dimensional moment problem. The technique suggested should
apply also to
other types of infinite-dimensional moment problems.
A similar approach was exploited recently in
the finite-dimensional (locally compact) case based on the the dual cone approach
in \cite{CuFi08}. We discuss this nonlocally compact case in
Section \ref{gensuff}, treating first the case of particle systems in
finite volume ($X$ compact) and give two alternative results in the
infinite volume ($X$ locally compact).

In Section \ref{limits} we derive mild conditions under which the
limit of
realizable correlation functions stays realizable. In
Section \ref{homogeneous} we show that correlation functions with some
symmetry, for example, translation invariance can, under mild extra
assumptions, be
realized by a point process with the same symmetry. In
Section \ref{secclass} we study a particular three-parameter subfamily of
the full set of necessary conditions derived earlier; we show that certain
well-known realizability conditions may be obtained from those of this
subfamily and that, in fact, they subsume all conditions from the
subfamily.

We now discuss briefly some applications of the realizability problem. As
already mentioned, the
problem has a long tradition in the theory of classical fluids
\cite{Percus,Percus2,Hansen}. It arises there because an
important ingredient of the theory is the introduction of various
approximation schemes, such as the Percus--Yevick and hyper-netted chain
approximations \cite{Hansen}, for computing the first two correlation
functions of the positions of the fluid molecules. It is then of interest
to determine whether or not the resulting functions in fact correspond to
any point process, that is, are in some sense internally consistent. If
they are, then they provide rigorous bounds for properties of the system
under consideration. The realizability problem was extensively discussed
in \cite{ST,ST2,ST4,ST5}
which consider the realization problem in various contexts, including a
conjecture related to the problem of the maximal density of sphere packing
in high dimensions \cite{ST3}.

The quantum mechanical variant of the realizability problem, known as
representability problem for reduced density matrices, is the basis of one
approach to the computation of the ground state energies of molecules
\cite{L,M,Cole1,Cole2,ColeYuk} yielding rigorous lower bounds. Interest
in this method is rising at present because improved algorithms in
semi-definite programming have led to an accuracy superior to that of the
traditional electronic structure method. These new methods are numerically
robust and reproduce further properties of the ground state; they are,
however, at present not competitive in terms of computation time
\cite{NewRepresent}. In \cite{GarrodPercus,Kummer}, the authors
give sufficient conditions for representability for systems with a fixed
finite number of particles, based on the dual cone approach mentioned
above. (Reference \cite{Percus2} gives corresponding classical results; see
Remark \ref{compsuff}(b) below.)

Applications of the problem of describing a point process from its low
order correlations also occur in biological contexts; for example,
in spatial ecology \cite{ecolo} and in the
study of neural spikes \cite{Mitra1,Mitra2}. In this and other situations
it is natural to consider a closely related problem in which the
correlation functions are specified only on part of the domain $X$; for
example, if $X$ is a lattice then we might only specify the nearest
neighbor correlations. See \cite{Kanter} for a similar problem in error
correcting codes. This will not be considered here; see, however,
\cite{KLS}, Section 7.

\section{Definitions}\label{ptproc}

We consider point processes in a
locally compact space $X$ which has a countable base of the topology.
$X$ is then a complete separable metric space under an appropriate metric
$d_X$ \cite{DV-J}, that is, it is a Polish space. We will sometimes use
the fact that such a metric exists for which closed balls of finite radius
are compact \cite{HY}. Measurability in $X$ will for us always mean
measurability with respect to the Borel $\sigma$-algebra on $X$. We will
write $\C_c(X)$ for the space of real-valued continuous functions with
compact support on~$X$ and $\M_c(X)$ for the space of real-valued bounded
measurable functions with compact support on $X$. The specific examples
that we have in mind for $X$ include the Euclidean spaces $\RR^d$,
manifolds (in particular the torus) and countable sets equipped with the
discrete topology. In the following we refer for brevity to these
countable sets as \textit{lattices}; the important special cases are $\ZZ^d$
and the discrete toruses. For the spaces $\RR^d$,~$\ZZ^d$ and the usual
and the discrete toruses one has as additional structure: a natural action
of the group of translations and the (uniform) measure which is invariant
under this action.

Intuitively, a point process on $X$ is a random distribution of points
in $X$ such that, with probability one, any compact set contains only
finitely many of these points. To give a precise definition, recall
that a \textit{Radon measure} on~$X$ is a Borel measure which is finite
on compact sets and denote by~$\mathcal{N}(X)$ the space of all Radon
measures $\eta$ on $X$ which take as values either a~nonnegative
integer (i.e., a member of $\NN_0=\{0,1,\ldots\}$) or infinity. A
measure $\eta\in\N(X)$ corresponds to a point configuration via the
representation
%
%
\begin{equation}
\label{eqdirac} \eta(d\mathbf{r}) = \sum_{i \in I} \delta_{\x
_i}(d\mathbf{r}),
\end{equation}
where either $I$ is finite or $I=\NN:= \{1,2,\ldots\}$; $\x_i \in X$ for
$i\in I$ and if $I=\NN$ the sequence $(\x_i)_{i\in I}$ has no accumulation
points in $X$; and $\delta_{\x_i}$ is the unit mass (Dirac measure)
supported at $\x_i$. Note that in this formulation there can be several
distinctly labeled points of the process at the same point of $X$. The
correspondence between $\eta$ and $(\x_i)_{i\in I}$ is one-to-one modulo
relabeling of the points. The requirement that $\eta$ be a Radon measure
corresponds to the condition that any compact set contain only finitely
many points of the process.

We equip $\mathcal{N}(X)$ with the vague topology which is the weakest
topology in which the mappings
%
%
\begin{equation}
\label{deflinf} \eta\mapsto\langle f, \eta\rangle:= \int_X f(\mathbf{r})
\eta(d\mathbf{r})
\end{equation}
are continuous for all $f \in\mathcal{C}_c(X)$. $\mathcal{N}(X)$ with
this topology is a Polish spa\-ce~\cite{DV-J}. Then we define a \textit{point
process} to be a Borel probability measure $\mu$ on~$\mathcal{N}(X)$. If
$\Ns$ is a measurable subset of $\N(X)$ $\mu(\mathcal{N}_{\supp})=1$,
we will
say that $\mu$ is a \textit{point process on $\Ns$}.

When $X$ is a lattice, $\mathcal{N}(X)$ can be identified with $\NN_0^X$
equipped with
the product topology; $\eta\in\N(X)$ is then identified with the function
on $X$ for which $\eta(\mathbf{r})$ is the number of particles at the
site $\mathbf{r}$. A
special case is the so-called lattice gas\label{latgas} in which there
can be at most
one particle per site, that is, $\eta(\mathbf{r}) \in\{0,1\}$. On the
lattice, of
course, integrals in formulas like (\ref{deflinf}) become sums, the Dirac
measure $\delta_\x(d\mathbf{r})$ becomes a Kronecker delta function,
etc. We will
not usually comment separately on the lattice case, adopting notation as
in (\ref{deflinf}) without further comment.

One advantage of defining point configurations as Radon measures is the
ease of then defining powers of these configurations. For $\eta\in\N(X)$,
$\eta^{\otimes n}$ denotes the (symmetric Radon) product measure on $X^n$;
note that from (\ref{eqdirac}) we have
%
%
\begin{equation}
\eta^{\otimes n}(d\mathbf{r}_1,\ldots,d\mathbf{r}_n)
= \sum_{i_1,i_2,\ldots,i_n}
\prod_{k=1}^n\delta_{x_{i_k}}(d\mathbf{r}_k).
\end{equation}
Here we will use a notation parallel to (\ref{deflinf}): for
$f_n\dvtx X^n \rightarrow\RR$ measurable and nonnegative, or for
$f_n\in\M_c(X^n)$, we write
%
%
\begin{eqnarray}\label{defprod2}
\langle f_n, \eta^{\otimes n} \rangle
:\!&=& \int_{X^n} f_n(\mathbf{r}_1,\ldots,\mathbf{r}_n)\eta(d\mathbf{r}_1)\cdots
\eta(d\mathbf{r}_n)\nonumber\\[-8pt]\\[-8pt]
&=& \sum_{i_1,i_2,\ldots,i_n}
f_n(\x_{i_1},\ldots, \x_{i_n}).\nonumber
\end{eqnarray}
By convention, $\langle f_0,\eta^{\otimes0}\rangle=f_0$ for $f_0\in
\RR$.
We will occasionally use a similar notation for functions: if $f\dvtx X\to
\RR$,
then $f^{\otimes n}(\mathbf{r}_1,\ldots,\mathbf{r}_n)=f(\mathbf{r}_1)\cdots
f(\mathbf{r}_n)$.

We will also need the \textit{factorial $n$th power} $\eta^{\odot n}$
of $\eta$, the symmetric Radon measure on $X^n$ given by
%
%
\begin{equation}\label{lstmin}
\eta^{\odot n}(d\mathbf{r}_1,\ldots,d\mathbf{r}_n)
:= \mathop{{\sum}'}_{i_1,i_2, \ldots, i_n}
\prod_{k=1}^n\delta_{x_{i_k}}(d\mathbf{r}_k),
\end{equation}
where $\sum'$ denotes a sum over distinct indices
$i_1,i_2, \ldots, i_n$, so that, in parallel to~(\ref{defprod2}),
%
%
\begin{eqnarray}\label{defnreod}
\langle f_n, \eta^{{\odot n}} \rangle
&=& \int_{X^n} f_n(\mathbf{r}_1,\ldots,\mathbf{r}_n)\eta^{\odot n}(d\mathbf{r}_1,\ldots,
d\mathbf{r}_n) \nonumber\\[-8pt]\\[-8pt]
&=& \mathop{{\sum}'}_{i_1, i_2, \ldots, i_n}
f_n(\x_{i_1},\ldots, \x_{i_n}).\nonumber
\end{eqnarray}
The term ``factorial power'' arises because, for any measurable subset
\mbox{$A$ of~$X$},
%
%
\begin{equation}\label{factorial}
\langle\mathbf1_A^{\otimes n},\eta^{\odot n}\rangle
\equiv\eta^{\odot n}(A\times\cdots\times A)
=\eta(A)\bigl(\eta(A)-1\bigr)\cdots\bigl(\eta(A)-n+1\bigr).
\end{equation}

One may view $\N(X)$ as a subset (with the inherited topology) of the
vector space of all signed Radon measures on $X$, equipped again with
the vague topology. Motivated by this imbedding we call functions on
$\N(X)$ of the form (\ref{deflinf}) \textit{linear}, since they are the
restrictions to $\N(X)$ of linear functionals. More generally, we
define a \textit{polynomial} on $\N(X)$ to be a function of the form
%
%
\begin{equation}
\label{defmonom} P(\eta) := \sum_{m=0}^n\langle f_m,\eta^{\odot
m}\rangle,
\end{equation}
where $f_0\in\RR$ and $f_m\in\M_c(X^m)$, $m=1,\ldots,n$; without loss
of generality we will assume that $f_m$ is symmetric in its arguments
when $m\ge2$. [We would obtain the same set of polynomial functions if
in (\ref{defmonom}) we replaced $\eta^{\odot m}$ by~$\eta^{\otimes
m}$.] We will sometimes consider \textit{polynomials with continuous
coefficients}, that is, polynomials for which $f_m\in\C_c(X)$,
$m=1,\ldots,n$.

\subsection{Correlation functions} \label{seccor}

It is often convenient to study point processes through their
\textit{correlation measures}, also called \textit{factorial moment
measures} or \textit{correlation functions}. The $n$th correlation
measure is the expectation of the $n$th factorial power:
%
%
\begin{equation}\label{defrhon0}
\rho_n(d\mathbf{r}_1,\ldots,d\mathbf{r}_n) \dvtx
= \mathbb{E}_\mu[ \eta^{\odot n}(d\mathbf{r}_1,\ldots,d\mathbf{r}_n) ],
\end{equation}
that is, it is the symmetric measure $\rho_n$ on $X^n$ satisfying
%
%
\begin{equation}\label{defrhon}
\int_{X^n} f_n(\mathbf{r}_1,\ldots,\mathbf{r}_n) \rho_n(d\mathbf{r}_1,\ldots
,d\mathbf{r}_n)
=\int_{\mathcal{N}(X)} \langle f_n, \eta^{\odot n} \rangle\mu
(d\eta)
\end{equation}
for all nonnegative measurable functions $f_n$ on $X^n$. One may also
define the $n$th \textit{moment measure} of the process by replacing
$\eta^{\odot n}$ by $\eta^{\otimes n}$ in (\ref{defrhon0}) and (\ref{defrhon})
but these measures will not play a significant role in our discussion. The
two sorts of moment measures are easily related; for example, at first
order they coincide, since $\eta^{\odot1}=\eta^{\otimes1}=\eta$
and at second
order we have
\begin{eqnarray*}
\int_{X\times X} f_2(\mathbf{r}_1,\mathbf{r}_2) \rho_2(d\mathbf{r}_1,d\mathbf{r}_2) &=& \int_{\mathcal{N}(X)}
\int_X \int_X f_2(\mathbf{r}_1,\mathbf{r}_2) \eta(d\mathbf{r}_1) \eta(d\mathbf{r}_2) \mu(d\eta) \\
&&{} - \int_{\mathcal{N}(X)} \int_X f_2(\mathbf{r},\mathbf{r})
\eta(d\mathbf{r}) \mu(d\eta).
\end{eqnarray*}

We will usually refer to the $\rho_n$ as \textit{correlation
functions} since this is the standard terminology in the physics
literature. This usage is particularly appropriate on a lattice or when
the measures are absolutely continuous with respect to Lebesgue
measure, if
we then gloss over the distinction between a~measure and its density.
From a more general viewpoint the terminology can be justified considering
$\rho_n$ as a generalized function in the sense of Schwartz. When $X$
is a
lattice, the process is a lattice gas; cf. page \pageref{latgas}, if
and only if
$\rho_2(\mathbf{r},\mathbf{r})=0$ for each $\mathbf{r}\in X$.

We say that the point process $\mu$ has \textit{finite local
$n$th moments} if for every compact subset $\Lambda$ of $X$,
%
%
\begin{eqnarray}\label{conmom}
\mathbb{E}_\mu[\eta(\Lambda)^n]
&\equiv&\mathbb{E}_\mu[\langle{\mathbf1}_\Lambda, \eta\rangle^n]
\equiv\mathbb{E}_\mu[\langle{\mathbf1}_\Lambda^{\otimes n},
\eta^{\otimes n} \rangle]\nonumber\\[-8pt]\\[-8pt]
&\equiv&\int_{\mathcal{N}(X)} \eta(\Lambda)^n \mu(d\eta) <
\infty.\nonumber
\end{eqnarray}
Obviously, the point process has then also finite local $m$th
moments for all $m\leq n$. If (\ref{conmom}) holds for $\Lambda=X$ we say
that $\mu$ has \textit{finite $n$th moment}. It is easy to see that
(\ref{conmom}) is equivalent to
$\rho_{n}(\Lambda^{n})\equiv\mathbb{E}_\mu[\langle
{\mathbf1}_\Lambda^{\otimes n}, \eta^{\odot
n}\rangle]<\infty$ [e.g., this follows by taking $\chi={\mathbf
1}_\Lambda$
in (\ref{ineq}) below]; in other words, the correlation measures $\rho_{m}$
are $\sigma$-finite Radon measures for all $m \leq n$ if and only if
(\ref{conmom}) holds. When the process has finite local $n$th moment
one may extend (\ref{defrhon}) to all $f_m\in\M_c(X^m)$ for $m \leq n$.
In this paper we will assume, unless it is specifically stated otherwise,
that the point processes under consideration have finite local second
moments.

\section{The realizability problem}\label{realize}

In Section \ref{seccor} we discussed how a point process $\mu$ gives rise
to correlation functions $\rho_n$. The realizability problem is a sort of
inverse problem.
\begin{definition} \label{realization}
Given $N \in\mathbb{N}$, symmetric Radon measures $\rho_n$ on $X^n$ for
$n=1,\ldots, N$ and a measurable subset $\Ns$ of $\mathcal{N}(X)$, we
say that $(\rho_n)_{n=1,\ldots,N}$ is \textit{realizable on $\Ns$} if
there exists a point process $\mu$ on $\Ns$ which for $n=1,\ldots,N$
has $\rho_n$ as its $n$th correlation function.
\end{definition}

Notice that, because the $\rho_n$ in Definition \ref{realization} are
assumed to be Radon measures, the realizing measure $\mu$ must have finite
local $N$th moments.

The aim of this paper is to develop necessary and sufficient conditions for
realizability solely in terms of $(\rho_n)_{n=1,\ldots,N}$. We will
describe these conditions in detail for the case $N=2$; the generalization
to general $N$ is straight forward. The case $N=\infty$ was treated in
\cite{Le75a,Le75b}; the problem with $N$ finite involves certain additional
difficulties, one of which is that the realizing measure is now
generically nonunique (see also Example \ref{funny} and Remark~\ref{lenard}).

\subsection{Necessary conditions}\label{subsecnec}

It is rather easy to give very general necessary conditions for the
realizability problem. Let $P(\eta)$ be a quadratic polynomial on $\N(X)$,
%
%
\begin{equation}\label{quadpoly} P(\eta)=P_{f_0,f_1,f_2}(\eta)
:=f_0+\langle f_1,\eta\rangle+\langle f_2,\eta^{\odot2}\rangle.
\end{equation}
Let $\mu$ be a point process on a given $\Ns\subset\N(X)$;
according to
(\ref{defrhon}) the expectation $\EE_\mu[P]$ can be computed in
terms of
the first two correlation functions of $\mu$ as
%
%
\begin{equation}\label{condadmis}
\EE_\mu[P_{f_0,f_1,f_2}] = f_0 + \int_X f_1(\mathbf{r}) \rho_1(d\mathbf{r}) +
\int_{X^2} f_2(\mathbf{r}_1,\mathbf{r}_2) \rho_2(d\mathbf{r}_1,d\mathbf{r}_2).
\end{equation}
On the other hand, if $P_{f_0,f_1,f_2}$ is nonnegative on $\Ns$, that
is, if for all $\eta=\sum_{i\in I}\delta_{\x_i}\in\Ns$,
%
%
\begin{equation}
f_0+\sum_i f_1(\x_i)+\sum_{i\ne j}f_2(\x_i,\x_j)\ge0,
\end{equation}
then necessarily $\EE_\mu[P_{f_0,f_1,f_2}]\ge0$. This leads
immediately to the following theorem.
\begin{theorem}[(Necessary conditions)] \label{lemnes}
If the pair $(\rho_1,\rho_2)$ is realizable by a point process on
$\Ns\subset\N(X)$ then for any quadratic polynomial $P_{f_0,f_1,f_2}$
which is nonnegative on $\Ns$,
%
%
\begin{equation}
\label{condness} f_0 + \int_X f_1(\mathbf{r}) \rho_1(d\mathbf{r}) +\int_{X^2}
f_2(\mathbf{r}_1,\mathbf{r}_2) \rho_2(d\mathbf{r}_1,d\mathbf{r}_2) \geq0.
\end{equation}
\end{theorem}

Theorem \ref{lemnes} gives uncountably many necessary conditions for
realizability indexed by the triples $(f_0,f_1,f_2)$. In
Section \ref{secclass} we will discuss how various standard conditions for
realizability are obtained from one class of such triples. Unfortunately,
the practical use of the theorem is limited because it is very
difficult to
identify admissible triples which lead to new and useful necessary
conditions.

\subsection{Sufficient conditions: Hard core exclusion}\label{susechard}

The idea of a ``hard core exclusion,'' which prevents the points of a
process from being too close together, is a common one in statistical
physics. To be precise:
\begin{definition}\label{hardcore}
Suppose that $d$ is a metric for the topology of $X$ and $D>0$. A
symmetric measure $\rho_2$ on $X\times X$ \textit{forces a hard core
exclusion with diameter $D$ for the metric $d$} if
%
%
\begin{equation}\label{hccond}
\rho_2\bigl(\{ (\mathbf{r}_1,\mathbf{r}_2)\in X\times X\mid d(\mathbf{r}_1,\mathbf{r}_2)< D \}\bigr)=0.
\end{equation}
\end{definition}

Condition (\ref{hccond}) says that, with probability one, no two
points of the process can lie in a distance less than $D$ from each other.
It is clear that if $\rho_2$ forces a hard core exclusion with
diameter $D$
then any point process with second correlation
function $\rho_2$ must be supported on
\[
\N_D(X) :=\biggl\{\eta\equiv\sum_i
\delta_{\x_i}\bigm| d (\x_i,\x_j)\ge D \mbox{ for all } i\neq
j\biggr\}.
\]
In this
subsection we show that under this hard core hypothesis the
necessary condition of Section \ref{subsecnec} for realizability on
$\N_D(X)$ is also sufficient.
\begin{theorem} \label{thsuflat}
Let $(\rho_1,\rho_2)$ be Radon measures on $X$ and $X\times X$,
respectively, with $\rho_2$ symmetric and suppose that $\rho_2$ forces
a hard core exclusion with diameter $D$ for a metric $d$. Then
$(\rho_1,\rho_2)$ is realizable on $\N_D(X)$ if and only if for any
quadratic polynomial $P_{f_0,f_1,f_2}(\eta)$ which is nonnegative on
$\N_D(X)$, $f_0$, $f_1$ and $f_2$ satisfy (\ref{condness}).
\end{theorem}
\begin{remark}\label{hcremark}
(a) The hard core exclusion
condition of
Definition \ref{hardcore} depends on the choice of metric $d$. Note,
however, that if $\rho_2$ satisfies (\ref{hccond}) for some metric
(generating the topology of $X$) then $(\rho_1,\rho_2)$ will be realizable
on the domain $\N_D(X)$ defined using that metric. In the following we will
not stress the dependence of $\N_D(X)$ on the metric.

(b) If $X$ is a lattice, then a point process realizing $(\rho_1,\rho_2)$
is a lattice gas if and only if there exists a metric $d$ and a $D>0$ such
that $\rho_2$ forces a hard core exclusion with diameter $D$ for the metric
$d$. If $\rho_2$ forces a hard core exclusion for some $d$ and $D$,
then certainly $\rho_2(\mathbf{r},\mathbf{r})=0$ for all $\mathbf{r}$; on the
other hand, given a
lattice gas we may topologize $X$ via the metric in which $d(\mathbf{r}_1,\mathbf{r}_2)=1$
whenever $\mathbf{r}_1\ne\mathbf{r}_2$ and in this metric $\rho_2$ forces
an exclusion
with diameter $D=1/2$. Thus, for lattice gases, Theorem \ref{thsuflat} gives
necessary and sufficient conditions for realizability with $\N_D(X)$ just
the set of lattice gas configurations. Of course, other hard core
restrictions are possible; for example, on $\ZZ^d$ in the standard metric
we may, in this way, forbid simultaneous occupancy of two nearest neighbor
sites.

(c) If $X$ is a finite set and $\rho_2$ forces an exclusion via
$\rho_2(\mathbf{r},\mathbf{r})\equiv0$ then~$\N(X)$ is finite and the
question of
realizability is one of the feasibility of a (finite) linear programming
problem: to find $(p_\eta)_{\eta\in\N(X)}$ with $p_\eta\ge0$ and, for
$\mathbf{r},\mathbf{r}_1,\mathbf{r}_2\in X$ with $\mathbf{r}_1\ne\mathbf{r}_2$,
\[
\sum_\eta p_\eta=1,\qquad
\sum_{\eta(\mathbf{r})=1} p_\eta=\rho_1(\mathbf{r}) \quad\mbox
{and}\quad
\sum_{\eta(\mathbf{r}_1)=\eta(\mathbf{r}_2)=1} p_\eta=\rho_2(\mathbf{r}_1,\mathbf{r}_2).
\]
By the duality theorem of linear programming the problem is feasible if
and only if a certain dual minimization problem has nonnegative solution.
But in fact the dual problem involves the coefficients of what we have
called a~quadratic polynomial in~$\eta$, the constraints of the problem
correspond to the positivity of this polynomial and the quantity to be
minimized is just the left-hand side of (\ref{condness}); that is,
Theorem~\ref{thsuflat} is equivalent in this case to the duality theorem.
The realization problem on a finite set can thus be studied numerically
via standard linear programming methods (see, e.g. \cite{CKLS}).
\end{remark}

For convenience we collect here three standard results which will be
used in proving Theorem \ref{thsuflat} and in Section \ref{gensuff}. In
stating the first two we will let $V$ be a vector space of real-valued
functions on a set $\Omega$. $V$ is a \textit{vector lattice} if for
every $v \in V$ also $|v| \in V$ (equivalently $v_+ \in V$). On $V$ we
may consider the natural (pointwise) partial order; we say that a
subspace $V_0$ of $V$ \textit{dominates} $V$ if for every $v \in V$
there exist $v_1,v_2 \in V_0$ such that $v_1 \leq v \leq v_2$. Then
\cite{Le75b,AhKr62}:
\begin{theorem}[(Riesz--Krein extension theorem)]\label{RKET}
Suppose that $V$ is a~vector space of functions as above and let $V_0$
be a subspace that dominates $V$. Then any nonnegative linear
functional on $V_0$ has at least one nonnegative linear extension to
all of $V$.
\end{theorem}

We note that the nonuniqueness of the extension given by this theorem
is the root of the nonuniqueness, mentioned above, of the realizing point
process. The next result is from the Daniell theory of integration
\cite{Ro62,P67}:
\begin{theorem}\label{theodaniel}
Let $V$ be a~vector space of functions as above which is a~vector
lattice and which contains the constant functions. Let $L$ be a~%
nonnega\-tive linear functional on $V$ for which:
\begin{enumerate}[(D)]
\item[(D)] If $(v_n)_{n \in\mathbb{N}}$ is a sequence of functions
in $V$
which decreases monotonically to zero then
$\lim_{n \rightarrow\infty} L(v_n)=0$.
\end{enumerate}
Then there exists
one and only one measure $\nu$ on $(\Omega,\Sigma_V)$, where $\Sigma
_V$ is
the $\sigma$-algebra generated by $V$, such that for all $v\in V$,
\[
L(v) = \int_\Omega v(\omega) \nu(d\omega).
\]
\end{theorem}

Finally we give a well-known characterization of compact subsets of
$\N(X)$ which follows from \cite{DV-J}, Corollary A.2.6.V and the
observation in Section \ref{ptproc} that $X$ is metrizable with a metric for
which all bounded sets have compact closure.
\begin{lemma} \label{lemcomp}
A set $C \subset\mathcal{N}(X)$ is compact if and only if $C$ is closed
and $\sup_{\eta\in C}\eta(\Lambda)< \infty$ for every compact subset
$\Lambda\subset X$.
\end{lemma}

Our next result is the key step in the proof of Theorem \ref{thsuflat}.
\begin{proposition} \label{propsuff} Let $\Ns$
be a compact subset of $\mathcal{N}(X)$, let $(\rho_1,\rho_2)$ be Radon
measures on $X$ and $X\times X$, respectively, with $\rho_2$ symmetric and
suppose that any quadratic polynomial $P_{f_0,f_1,f_2}(\eta)$ which is
nonnegative\break on~$\Ns$ satisfies~(\ref{condness}). Then $(\rho
_1,\rho_2)$
is realizable by a point process supported on
$\mathcal{N}_{\mathrm{supp}}$.
\end{proposition}
\begin{pf}
Let $V$ be the vector space of all continuous functions on
$\Ns$ and let $V_0$ be the vector space of all quadratic polynomials
$P_{f_0,f_1,f_2}$ with continuous coefficients; from the compactness of
$\Ns$ it is clear that $V_0$ dominates $V$. Let $L$ be the linear form on
$V_0$ defined by
\[
L(P_{f_0,f_1,f_2}) := f_0 + \int_X f_1(\mathbf{r}) \rho
_1(d\mathbf{r}) +
\int_{X \times X} f_2(\mathbf{r}_1,\mathbf{r}_2) \rho_2(d\mathbf{r}_1,d\mathbf{r}_2).
\]
The hypothesis of the theorem is precisely that $L$ is nonnegative
so by the
Riesz--Krein extension theorem we can extend $L$ to a nonnegative linear
functional on all of $V$.
Since $\Ns$ is compact, the Riesz--Markov representation theorem implies
that there exists a probability
measure $\mu$ on $\Ns$---that is, a point process on $X$---such that
\[
L(F) = \int_{\mathcal{N}_{\mathrm{supp}}} F(\eta) \mu(d\eta)
\]
for all $F\in V$. In particular, taking
$F_n(\eta)=\langle f_n,\eta^{\odot n}\rangle$ for $n=1,2$, with
$f_n\in\C_c(X^n)$ and $f_2$ symmetric, we obtain (\ref{defrhon}) for
$n=1,2$ for continuous $f_1,f_2$; this suffices to imply that $\rho_1$ and
$\rho_2$ are indeed the correlation functions of the process~$\mu$.
\end{pf}

Note that the proof shows that it suffices for realizability that
(\ref{condness}) holds for polynomials with continuous coefficients.
\begin{pf*}{Proof of Theorem \ref{thsuflat}}
If $\mu$ is a realization of
$\rho_1,\rho_2$ then, as observed above, it must be supported on $\N_D(X)$,
and by Theorem \ref{lemnes} must satisfy the given condition. As the set
$\mathcal{N}_D(X)$ is compact, by Lemma \ref{lemcomp}, the
converse direction follows from Proposition \ref{propsuff}.
\end{pf*}

Hard core exclusion is not the only natural possibility for a compact
$\Ns$. If $N$ is a natural number then the set of all configurations
with exactly $N$ particles, or at most $N$ particles,
\begin{eqnarray*}
\N^N(X)
&:=&\{ \eta\in\N(X)\mid \eta(X)=N \}
,\\
\N^{\leq N}(X)
&:=&\{ \eta\in\N(X)\mid \eta(X)\leq N \},
\end{eqnarray*}
is compact. We summarize the consequences in the following corollary.
\begin{corollary}\label{CoNrep} Let $(\rho_1,\rho_2)$ be Radon
measures on
$X$ and $X\times X$ with $\rho_2$ symmetric. Suppose that
any quadratic polynomial $P_{f_0,f_1,f_2}(\eta)$ which is nonnegative on
$\N^N(X)$ [resp., $\N^{\leq N}(X)$] satisfies (\ref
{condness}). Then
$(\rho_1,\rho_2)$ is realizable by a point process supported on $\N^N(X)$
[resp., $\N^{\leq N}(X)$].
\end{corollary}

A similar result would hold for $X$ a lattice and, for some $k\ge0$,
$\Ns$
the set of configurations with at most $k$ particles at any site.
\begin{remark}\label{compsuff}
(a) The essential property for the proof
of Proposition~\ref{propsuff} is the compactness of $\Ns$. Indeed, the
result is false if $\Ns$ is replaced by $\N(X)$; see Example \ref{funny}
below.

(b) Corollary \ref{CoNrep} was established by Percus in \cite{Percus} and
\cite{Percus2} using the technique of double dual cone. This technique
should give an alternative approach to prove sufficiency of the conditions
but will require a careful identification of the closure of the initial
cone requiring considerations similar to those above. In
\cite{GarrodPercus} and \cite{Kummer} a quantum mechanical version of
Corollary \ref{CoNrep} was worked out in the framework of reduced density
matrices and trace class operators. A characterization of the closure of
the cone was not considered.

(c) For any given $(\rho_1,\rho_2)$ one could, of course, attempt to use
Proposition \ref{propsuff} to establish realizability on some suitably chosen
compact subset $\Ns\subset\N(X)$. For translation invariant
$(\rho_1,\rho_2)$ in $\RR^d$, for example, one might require that for
$\Lambda\subset X$ with volume $|\Lambda|$,
$\eta(\Lambda)\le A(1+|\Lambda|^k)$ for suitably chosen $A$ and $k$.
We do
not, however, know of an example in which such an approach succeeds. What
is significant about processes with hard cores is that the hard core
constraint is of physical interest, is expressible in terms of the given
datum $\rho_2$ and forces any realization to be on a compact set of
configurations.
\end{remark}

\subsection{Sufficient conditions without a hard core}\label{gensuff}

We now consider the case of general $(\rho_1,\rho_2)$ in which we
have no
{a priori} reason, such as a hard core constraint, to expect a
realizing process to be supported on a compact set of configurations. In
this case the necessary conditions of Theorem \ref{lemnes} are in general
not sufficient, as shown by the following example.
\begin{example}\label{funny}
Let $X=\RR^d$ and consider the pair of correlation functions
$\rho_1(\mathbf{r})\equiv0$,
$\rho_2(\mathbf{r}_1,\mathbf{r}_2)\equiv1$. This is certainly not
realizable, since if it were realized by some process $\mu$ then for
any measurable set $\Lambda$,
$\EE_\mu[\eta(\Lambda)]=\int_\Lambda\rho_1(d\mathbf{r})=0$ and hence,
$\eta(\Lambda)=0$ with probability one so that the second correlation
function of $\mu$ would have to vanish. But consider the point process
%
%
\begin{equation}
\mu^{\epsilon}(d\eta)
:= (1-\epsilon^2) \delta_0(d\eta) + \epsilon^2 \pi_{1/\epsilon
}(d\eta),
\end{equation}
where $\epsilon\in(0,1]$, $\pi_z$ denotes the Poisson measure on
$\RR^d$
with density $z$ and $\delta_0$ is the measure concentrated on $\eta=0$.
The corresponding correlation functions
$\rho_1^{\varepsilon}(\mathbf{r}) = \varepsilon$ and
$\rho_2^{\varepsilon}(\mathbf{r}_1,\mathbf{r}_2) = 1$ converge as $\epsilon
\to0$ to the
given $(\rho_1,\rho_2)$, from which it follows easily that the latter
fulfills the necessary condition of Theorem \ref{lemnes}.
\end{example}

In the following subsections we give sufficient conditions for
realizability in the general case. Lemma \ref{lemcomp} and
Proposition \ref{propsuff} indicate that difficulties in doing so will be
associated with the local occurrence of an unbounded number of particles.
The key idea is to control this by requiring not only realization of
$\rho_1$ and $\rho_2$ but also the existence in some form of a finite third
moment (a moment of order $2+\epsilon$ would suffice). Such a requirement
can be motivated by reconsidering the proof of Theorem \ref{thsuflat},
omitting the hard core hypothesis and trying to prove existence of a
process supported on $\N(X)$. Defining $V$ to include only functions of
quadratic growth in $\eta$ and using Theorem \ref{theodaniel} rather than
the Riesz--Markov theorem, one may establish the existence of a process
$\mu$ realizing $\rho_1$ but not necessarily $\rho_2$. The situation in
this section (see, e.g., Theorem \ref{propsuffg}) is similar: by
controlling a third moment we can realize the first two correlation
functions. The condition can also be motivated by considering
Example \ref{funny}; no similar example can be constructed in which the
third moments of the processes $\mu^\epsilon$ are uniformly bounded.
\begin{remark}\label{lenard}
(a) Even if $X$ is a lattice one will still
need to control~so\-me higher moment if there is no bound on the number of
particles per site.

(b) In the case in which all correlation
functions are prescribed, that is, when $N=\infty$ in the sense of
Definition \ref{realization}, the need to control an ``extra'' moment does
not arise. See \cite{Le75a,Le75b}.
\end{remark}

Since the essential difficulties are local they will occur even for
compact $X$; we will first discuss this case where certain technical
difficulties are absent. Throughout this section we will
define the function $H^\chi_n$ on $\N(X)$, where $\chi$ is a
strictly positive
bounded continuous function on $X$ and $n\ge0$, by
%
%
\begin{equation}\label{defHchi}
H^\chi_n(\eta):=\langle\chi^{\otimes n},\eta^{\odot n}\rangle
= \mathop{{\sum}'}_{i_1,\ldots, i_n}\chi(x_{i_1})\cdots\chi(x_{i_n})
\end{equation}
with $\sum'$ as in (\ref{lstmin}). Note that since all summands in
(\ref{defHchi}) are nonnegative the sum always is well defined,
though it
may be infinite.
For $\Lambda\subset X$ we write
$H_n^\Lambda:= H_n^{{\mathbf1}_\Lambda}$ and we abbreviate $H^X_n$ as
$H_n$. In (\ref{defHchi}) we have defined $H^\chi_n$ using the factorial
power $\eta^{\odot n}$ but one could equivalently work with
$\eta^{\otimes n}$; this follows from the fact that for each $n\ge0$ there
exists a constant $b_n>0$ such that for all $\eta\in\N(X)$,
%
%
\begin{equation}\label{ineq}
\tfrac12 \langle\chi,\eta\rangle^n -b_n
\leq H_n^\chi(\eta) \leq\langle\chi,\eta\rangle^n
\equiv\langle\chi^{\otimes n},\eta^{\otimes n}\rangle.
\end{equation}
To verify (\ref{ineq}) we note that as all summands in (\ref
{defHchi}) are
nonnegative the upper bound is immediate. On the other hand, the
difference of $\langle\chi,\eta\rangle^n-H_n^\chi(\eta)$
can be bounded by a linear combination of $\langle \chi,
\eta \rangle^m$ for $m < n$ and each of
these can be estimated above by $c \langle\chi, \eta\rangle^n + c'$ for
$c>0$ arbitrary small.
As mentioned just below (\ref{conmom}), the\vspace*{1pt} inequalities (\ref{ineq}) implies
that $\mu$ has finite local
$n$th moments is equivalent to $\mathbb E_\mu[H^\Lambda
_n]<\infty$
for all compact $\Lambda$, so that $\mu$ has finite $n$th
moment if
and only if $\mathbb E_\mu[H_n]<\infty$.

We will say that $\mu$ has \textit{finite $n$th $\chi$-moment} if
$\mathbb E_\mu[H^\chi_n]<\infty$; in particular, $\mu$ then has
support on
the set of all configurations $\eta$ with
$\langle\chi,\eta\rangle< \infty$. By~(\ref{ineq}) and the
positivity of
$\chi$, finite $n$th $\chi$-moment implies finite local $n$th
moments. Clearly the converse will not hold for general $\chi$ but we
will show in Lemma~\ref{lemchimoment} below that a measure with finite
local $n$th moments has finite $n$th $\chi$-moment for an
appropriately chosen $\chi$.\vspace*{-2pt}

\subsubsection{Compact $X$} \label{compact}

Suppose that $X$ is compact. In the next theorem we give a condition which
is both necessary and sufficient for $(\rho_1,\rho_2)$ to be
realizable by a
process with a finite third moment. As a corollary we obtain a
sufficient condition for realizability of $(\rho_1,\rho_2)$. The conditions
that we will give involve cubic polynomials of a special form that we will
call \textit{restricted}. These have the form
%
%
\begin{equation}\label{eqadmnew}
Q_{f_0,f_1,f_2,f_3}(\eta)
= f_0 +\langle f_1,\eta\rangle+ \langle f_2,\eta^{\odot2}\rangle
+ f_3 H_3(\eta),
\end{equation}
where $f_0,f_3\in\RR$, $f_1\in\C_c(X)$ and $f_2\in\C_c(X^2)$ with $f_2$
symmetric.
\begin{theorem} \label{propsuffg} Let $X$ be compact. Then symmetric Radon
measures $\rho_1$ and $\rho_2$ on $X$ and $X\times X$ are realizable
by a
point process with a finite third moment if and only if there exists a
constant $R>0$ such that any restricted cubic polynomial
$Q_{f_0,f_1,f_2,f_3}$ which is nonnegative on $\N(X)$ satisfies
%
%
\begin{equation}
f_0 + \int_X f_1(\x) \rho_1(d\x) + \int_{X \times X} f_2(\x
,\y)
\rho_2(d\x,d\y) +f_3 R \geq0.
\end{equation}
\end{theorem}

We now have:\vspace*{-2pt}
\begin{corollary} \label{corrsuff} If the condition of
Theorem \ref{propsuffg} holds then the pair $(\rho_1,\rho_2)$ is
realizable.\vspace*{-2pt}
\end{corollary}
\begin{pf*}{Proof of Theorem \ref{propsuffg}} Let $V$ be the vector
space of all
continuous functions $F$ on $\mathcal{N}(X)$ such that $|F|\le C(1+H_3)$
for some constant $C>0$ and let $V_0$ be the subspace of $V$ consisting of
all restricted cubic polynomials. For any $R\ge0$ we may define a linear
functional $L_R$ on $V_0$ by
\begin{eqnarray*}
L_R(Q_{f_0,f_1,f_2,f_3})
&:=& f_0 + \int_X f_1(\mathbf{r})
\rho_1(d\mathbf{r}) \\[-2pt]
&&{}+ \int_{X \times X} f_2(\mathbf{r}_1,\mathbf{r}_2)
\rho_2(d\mathbf{r}_1,d\mathbf{r}_2)+ f_3 R.
\end{eqnarray*}
Then we must show that $\rho_1,\rho_2$ is realizable by a measure
with a finite third moment if and only if $L_R$ is nonnegative for some
$R>0$.

The condition is clearly necessary since if $\mu$ is such a realizing
measure and $Q_{f_0,f_1,f_2,f_3}\ge0$ then
%
%
\begin{equation}
L_{\mu(H_3)}(Q_{f_0,f_1,f_2,f_3})
= \int Q_{f_0,f_1,f_2,f_3}(\eta) \mu(d\eta)
\ge0.
\end{equation}
Suppose conversely then that $R$ is such that $L_R$ is nonnegative on $V_0$.
It is easily seen that $V_0$ dominates $V$, so that, by Theorem \ref{RKET},
$L_R$ has a~nonnegative extension, which we will also call $L_R$, to
all of
$V$. It remains to show that this extended linear form is actually given
by a measure.

Let $W$ be the subspace of $V$ consisting of those functions $F\in
V$ such that $|F|\le C(1+H_2)$ for some $C>0$. $W$ is a lattice
which generates the $\sigma$-algebra corresponding to the vague
topology because it contains all functions of the form $\langle f ,
\cdot\rangle$ with $f$ continuous. We wish to apply
Theorem \ref{theodaniel} to $L_R$ on $W$ and so must verify that
$L_R$ satisfies (D). Let $(F_n)_{n\in\mathbb{N}}$ be a
monotonically decreasing sequence in $W$ which converges pointwise
to $0$ and let $\epsilon>0$ be given. The sets $K_n:=\{\eta\in
\mathcal{N}(X)\mid F_n(\eta) \geq\epsilon[1+H_3(\eta)]\}$ are
closed because~$F_n$ and~$H_3$ are continuous. Moreover, $K_n$ is
compact because, since $F_n\in W$, $K_n$ is for some $C>0$ a subset
of $\{\eta\in\mathcal{N}(X)\mid C[1+H_2(\eta)] \geq\epsilon
[(1+H_3(\eta)]\}$, and the latter set is compact by
Lemma \ref{lemcomp} since $\eta(X)$ is bounded on it. Because the
$F_n$ decrease to zero pointwise, $\bigcap_nK_n=\varnothing$, so there
must exist an $N \in\mathbb{N}$ with $K_n=\varnothing$ for $n\ge N$,
that is, with $F_n \leq\epsilon(1+ H_3)$ for all $n \geq N$. This,
with the positivity of $L_R$, implies that for $n\geq N$,
\[
L_R(F_n) \leq L_R(Q_{\epsilon,0,0,\epsilon}) = \epsilon(1+R).
\]
As $\epsilon$ was arbitrary, (D) holds and, therefore,
Theorem \ref{theodaniel} implies that there exists a probability measure
$\mu$ on $\mathcal{N}(X)$ such that for $F\in W$,
\[
L_R(F) = \int_{\mathcal{N}(X)} F(\eta) \mu(d\eta).
\]
In particular, for all $f_0 \in\mathbb{R}$ and continuous
functions $f_1$ and $f_2$ on $X$ and $X\times X$,
\begin{eqnarray*}
&&
f_0 + \int_X f_1(\x) \rho_1(d\x) + \int_{X \times X} f_2(\x
,\y)
\rho_2(d\x,d\y) \\
&&\qquad = \int_{\mathcal{N}(X)}
P_{f_0,f_1,f_2}(\eta) \mu(d\eta),
\end{eqnarray*}
which implies that $\mu$ realizes $(\rho_1,\rho_2)$.

Finally, if for $n\in\NN$ we define $H^{(n)}_3(\eta)=\min\{H_3(\eta
),n\}$ then $H^{(n)}_3\in W$~and so $\int H^{(n)}_3(\eta) \,d\mu(\eta)
= L_R(H^{(n)}_3)$. But by the positivity of $L_R$ on $V$,
$L_R(H^{(n)}_3) \le\break L_R(H_3)= L_R(Q_{0,0,0,1})= R$ and so the monotone
convergence theorem implies that $\int H_3 \,d\mu\le R$, that is, $\mu$
has finite third moment.
\end{pf*}

\subsubsection{Noncompact $X$}\label{subsubnoncompact}

For the case in which $X$ is not compact we give, in Theorems
\ref{chisuff} and \ref{propnc}, two distinct sufficient conditions for
realizability which generalize Theorem \ref{propsuffg} in two different
ways. In this section we will finally assume that the metric $d_X$ is
such that bounded sets have compact closure; cf. the beginning of
Section \ref{ptproc}. With some fixed $\x_0\in X$ define $\LN=\{\x\in
X\mid d_X(\x,\x_0)\le N\}$. Throughout this section we let $\chi$ be a
strictly positive continuous function on~$X$. One should think of
$\chi$ as a function which vanishes at infinity; for example, if
$X=\RR^d$ we might take $\chi(\x)=(1+|\x|^2)^{-k}$ for some $k>0$.
\begin{lemma}\label{lemchimoment} A point process $\mu$ on $X$ has finite
local $n$th moments if and only if there exists a positive
continuous $\chi$ such that $\mu$ has finite $n$th $\chi$-moment.
\end{lemma}
\begin{pf} If $\mu$ has finite $n$th $\chi$-moment then, using
the continuity and positivity of $\chi$, it follows immediately that
$\mu$
has finite local $n$th moments. Suppose conversely that $\mu$ has
finite local $n$th moments. Let $\chi_k$ be a nonnegative function
on $X$ with compact support satisfying
$\Ii_{\Lambda_{k}}\le\chi_k \leq\Ii_X$; then
$\int_{\N(X)} \langle\chi_k , \eta\rangle^n \mu(d\eta) < \infty
$ for all
$k$. Define
\[
\chi(x) := \sum_{k=1}^\infty c_k \chi_k(x) \qquad\mbox{with }
c_k :=
\frac{2^{-k}}{1+\sqrt[n]{\int\langle\chi_k , \xi\rangle^n
\mu(d\xi)}}.
\]
Then
\begin{eqnarray*}
\int_{\mathcal{N}(X)} H_n^\chi(\eta) \mu(d\eta)
&\le&\int_{\N(X)}
\langle\chi, \eta\rangle^n \mu(d\eta)\\
&=& \sum_{k_1,k_2,
\ldots,k_n=1}^\infty\prod_{i=1}^nc_{k_i}
\int_{\N(X)}\prod_{i=1}^n\langle\chi_{k_i},\eta\rangle\mu
(d\eta) \\
&\leq&
\Biggl( \sum_{k=1}^\infty c_k \sqrt[n]{\int_{\N(X)} \langle\chi_{k}
, \eta\rangle^n \mu(d\eta)} \Biggr)^n \leq1,
\end{eqnarray*}
where we have used H\"older's inequality.
\end{pf}

The sufficiency criteria of the next theorem are stated in terms of
\textit{$\chi$-re\-stricted} cubic polynomials,
%
%
\begin{equation}\label{chirstdef}
Q^\chi_{f_0,f_1,f_2,f_3}(\eta)
= f_0 +\langle f_1,\eta\rangle+ \langle f_2,\eta^{\odot2}\rangle
+ f_3 H^\chi_3(\eta),
\end{equation}
where $f_0,\ldots,f_3$ are as in (\ref{eqadmnew}).
\begin{theorem} \label{chisuff} Symmetric Radon
measures $\rho_1$ and $\rho_2$ on $X$ and \mbox{$X\times X$} are realizable
by a
point process with finite local third moments if and only if there
exists a
constant $R>0$ and a positive function $\chi$ such that any
$\chi$-restricted cubic polynomial $Q^\chi_{f_0,f_1,f_2,f_3}$ which is
nonnegative on $\N(X)$ satisfies
%
%
\begin{equation}\label{chicond}
f_0 + \int_X f_1(\x) \rho_1(d\x) + \int_{X \times X} f_2(\x
,\y)
\rho_2(d\x,d\y) +f_3 R \geq0.
\end{equation}
\end{theorem}
\begin{pf} According to Lemma \ref{lemchimoment} it suffices to show
that $\rho_1$ and $\rho_2$ are realizable by a point process with finite
third $\chi$-moment if and only if the condition is satisfied. The
proof of
this is very similar to that of Theo\-rem~\ref{propsuffg}, with $H_2$ and
$H_3$ replaced by $H^\chi_2$ and $H^\chi_3$ throughout, so we content
ourselves here with commenting on the technical modifications necessitated
by the noncompact character of $X$.

One source of difficulties is that $H^\chi_n$ is not a
continuous function on $\N(X)$. This means that if the vector space used
in the proof was to be defined in parallel to the $V$ of the earlier proof
then it would not contain all $\chi$-restricted polynomials. The problem
may be avoided by replacing $V$ throughout by $V^\chi:=V^\chi
_0+V^\chi_1$,
where $V^\chi_0$ is the space of all $\chi$-restricted cubic
polynomials (which
plays the same role as did $V_0$ earlier) and $V^\chi_1$, defined in
parallel to the earlier $V$, is the vector space of all continuous
functions $F$ on $\mathcal{N}(X)$ such that $|F|\le C(1+H^\chi_3)$
for some
constant $C>0$.

The set $K_n$ is replaced by
$K^\chi_n=\{\eta
\in\mathcal{N}(X)\mid F_n(\eta) \geq\epsilon[(1+H^\chi_3(\eta)]\}
$; the argument that $K_n$ was closed used the
continuity of $H^\chi_3$ but lower semi-continuity suffices and we
establish this in the next lemma. $K^\chi_n$ is for some $C>0$ a
subset of
$\{\eta
\in\mathcal{N}(X)\mid C[1+H^\chi_2(\eta)] \geq\epsilon[1+H^\chi
_3(\eta)]\}$
and this set is precompact by Lem\-ma~\ref{lemcomp}, since $H^\chi_1$ is
bounded on it and for any compact $\Lambda\subset X$ there is a constant
$c_\Lambda$ with $\Ii_\Lambda\le c_\Lambda\chi$. The sequence $H^{(n)}_3$
used in the last step of the proof is replaced by any sequence of bounded
continuous functions increasing to $H^\chi_3$; the existence of such a
sequence follows from the lower semicontinuity of $H_3^\chi$.\vspace*{-2pt}
\end{pf}
\begin{lemma}\label{lsc} For any $n>0$ the function $H^\chi_n$ is lower
semi-continuous.\vspace*{-2pt}
\end{lemma}
\begin{pf} We must show that sets of the form
\mbox{$S\!:=\!\{ \eta\in\mathcal{N}(X) \mid H_n^\chi(\eta) \!\leq\! C\} $} are closed.
Let $(\eta_k)$ be a sequence in $S$ converging vaguely to
$\eta\in\mathcal{N}(X) $ and let $(\chi_m)$ be an increasing
sequence of
nonnegative continuous functions with compact support on $X$ such that
$\chi_m \nearrow\chi$. By the vague convergence
$\langle\chi_m^{\otimes n},\eta_k^{\odot n}\rangle\rightarrow
\langle\chi_m^{\otimes n},\eta^{\odot n}\rangle$
as $k\nearrow\infty$, for any fixed $m$, and
by the monotone convergence of $\chi_m$ also
$\langle\chi_m^{\otimes n},\eta^{\odot n}\rangle\nearrow H^\chi
_n(\eta)$
as $m\nearrow\infty$. Since $\langle\chi_m^{\otimes n},\eta
_k^{\odot
n}\rangle\leq H^\chi_n(\eta) \leq C$, also $H^\chi_n(\eta) \leq
C$, and so $S$ is closed.\vspace*{-2pt}
\end{pf}
\begin{remark}
By taking $\chi(\x)\equiv1$ in Theorem \ref {chisuff} we see that in
fact Theorem \ref{propsuffg} holds even when $X$ is not compact. For
typical problems, however, this result is not very interesting since a
realizing measure with finite third moment would be impossible if, for
example, $\langle\eta(X)\rangle=\int_X\rho_1(dx)$ were infinite, as
would be true for any nonzero translation invariant $\rho_1$ in~$\RR^d$.\vspace*{-2pt}
\end{remark}

We now give the second sufficient condition.\vspace*{-2pt}
\begin{theorem} \label{propnc} Let $X=\RR^d$. Then symmetric Radon
measures $\rho_1$ and~$\rho_2$ on $X$ and $X\times X$ are realizable
by a
point process with finite local third moments if and only if the condition
of Theorem \ref{propsuffg} holds in every $\Lambda_N$, $N\in\NN$.\vspace*{-2pt}
\end{theorem}
\begin{pf}
The necessity of the condition follows as in the proof of
Theorem~\ref{propsuffg}. Suppose conversely that the condition of
Theorem \ref{propsuffg} holds in every $\Lambda_N$ so that for each $N$
there exists a measure $\mu_N$ on $\N(\LN)$ which realizes
$(\rho_1,\rho_2)$ in $\LN$. If $N\ge n$ then $\mu_N$ defines in the
obvious way a marginal measure $\mu^n_N$ on $\N(\iLn)$, where
$\iLn$ denotes the interior of ${\Lambda_n}$; all the measures~$\mu^n_N$,
$N\ge n$, have the same one- and two-point correlation functions $\rho_1$
and $\rho_2$ on $\iLn$. Since
%
%
\begin{equation}
c_n:=\langle\eta(\iLn)\rangle_{\mu^n_N}=\int_{\iLn}\rho_1(dx)
\end{equation}
is independent of $N$, Markov's inequality implies that these measures
satisfy $\mu^n_N[(K_n(M)]\geq1- c_n/M$, where
$K_n(M)=\{\eta\mid\eta(\iLn)\le M\}$. Since $K_n(M)$ is compact by
Lemma \ref{lemcomp}, the sequence of measures $(\mu^n_N)_{N\ge n}$ is
tight and any subsequence of this sequence itself contains a convergent
subsequence. We may thus obtain recursively sequences
$(N_{n,k})_{k\in\NN}$ such that $(N_{n+1,k})$ is a subsequence of
$(N_{n,k})$ and such that $(\mu^n_{N_{n,k}})_{k\in\NN}$ converges
weakly to
a~measure $\mu^n$ on $\N(\iLn)$. The measure $\mu^n$ realizes
$(\rho_1,\rho_2)$ on $\iLn$. The $\mu^n$ are compatible, in the
sense that
$\mu^n$ is the marginal of $\mu^{n+1}$ on $\N(\iLn)$, because the
projections from $\N(\iLnp)$ onto $\N(\iLn)$ are continuous since
$\iLn$ is
open. Thus a realizing measure on $\N(X)$ exists by Kolmogorov's
projective limit theorem.
\end{pf}
\begin{remark}
In checking the sufficient conditions for realizability given in
Theorems \ref{thsuflat}, \ref{propsuffg}, \ref{chisuff} and
\ref{propnc} it may be advantageous to choose the coefficients of the
quadratic polynomials (\ref{quadpoly}) from a class of functions other
than $\C_c(X^m)$. Suppose then that we can verify the conditions
(\ref{condadmis}) when the coefficients $f_m$ are chosen from $\F_m$, a
subspace of $\C_c(X^m)$ [with $\F_0\equiv\C_c(X^0) \equiv\RR$]. By
straightforward modifications of the proofs of Theo\-rems \ref{thsuflat},
\ref{propsuffg} and \ref{chisuff} one sees that this will suffice for
realizability if $\F_m$ identifies measures on $X^m$, that is, if
whenever Radon measures $\nu$ and $\nu'$ satisfy $\int_{X^m}
f_m(\mathbf{r}) \nu(d\mathbf{r})= \int_{X^m} f_m(\mathbf{r}) \nu
'(d\mathbf{r})$ for all $f_m \in\F_m$, necessarily $\nu=\nu'$. For an
analogously modified version of Theorem \ref{propnc} slightly more is
needed: for each $N$ the functions from $\F_m$ with
$\operatorname{supp}{\F_m}\subset\LN$ must identify measures on $\LN^m$.
These conditions are fulfilled if $\F_m$ forms an algebra which
separates points. For example, if $X$ is a manifold without boundary
then one may take $\F_m$ to be $\C^\infty_c(X^m)$, the space of
infinitely differentiable functions with compact support.\vspace*{-3pt}
\end{remark}

\subsection{Stability of realizability under limits}\label{limits}

The sufficient conditions obtained above can be used to derive general
results about realizing measures. In this subsection we discuss sufficient
conditions for the limit of a sequence of realizable correlation functions
to be itself realizable. Each of the Theorems \ref{thsuflat},
\ref{propsuffg}, \ref{chisuff} and \ref{propnc}
will give rise to a different variant.
Recall that\vadjust{\eject} $(\rho_1^{(n)},\rho_2^{(n)})$
converges in the vague topology to $(\rho_1,\rho_2)$ if for any
$f_1\in\C_c(X)$ and $f_2\in\C_c(X\times X)$,\vspace*{-2pt}
%
%
\begin{equation}
\int_X f_1(\mathbf{r})\rho_1^{(n)}(d\mathbf{r})\to\int_X f_1(\mathbf{r})
\rho_1(d\mathbf{r})\vspace*{-2pt}
\end{equation}
and\vspace*{-2pt}
\begin{equation}
\int_{X2} f_2(\mathbf{r}_1,\mathbf{r}_2)
\rho_2^{(n)}(d\mathbf{r}_1,\mathbf{r}_2)\to\int_X
f_2(\mathbf{r}_1,\mathbf{r}_2) \rho_2(d\mathbf{r}_1,d\mathbf{r}_2).\vspace*{-2pt}
\end{equation}

For the hard core case we have to require a uniform exclusion diameter.\vspace*{-2pt}
\begin{proposition} \label{realcloscon} Let $(\rho_1^{(n)},\rho_2^{(n)})$
be a sequence of realizable pairs of symmetric Radon measures which
converges in the vague topology to $(\rho_1,\rho_2)$ and for which there
exists a $D>0$ such that
\mbox{$\rho_2^{(n)}(\{
(\mathbf{r}_1,\mathbf{r}_2)\mid d(\mathbf{r}_1,\mathbf{r}_2) < D \})=0$}
for all $n$. Then
$(\rho_1,\rho_2)$ is
also realizable.\vspace*{-2pt}
\end{proposition}
\begin{pf}
If
$P_{f_0,f_1,f_2}$ is a nonnegative quadratic polynomial on $\N_D(X)$
then the hypotheses imply that
%
%
\begin{equation}
f_0 + \int_X f_1(\x) \rho^{(n)}_1(d\x) + \int_{X \times X}
f_2(\x,\y) \rho^{(n)}_2(d\x,d\y) \geq0
\end{equation}
for all $n$. Taking the $n\to\infty$ limit then gives (\ref{condness}).
By the portmanteau~theo\-rem, the limiting correlation functions force
also a hard core exclusion.\vspace*{-2pt}
\end{pf}

For lattice gases this implies a very natural result:\vspace*{-2pt}
\begin{corollary} \label{realcloslat} Let $X$ be a lattice and let
$(\rho_1^{(n)},\rho_2^{(n)})$ be a sequence of realizable pairs with
$\rho_2^{(n)}(\mathbf{r},\mathbf{r})=0$ for all $n$ and $\mathbf{r}$. If
$(\rho_1^{(n)},\rho_2^{(n)})$ converges pointwise to $(\rho_1,\rho
_2)$, then
$(\rho_1,\rho_2)$ is realizable.\vspace*{-2pt}
\end{corollary}

From Theorem \ref{propsuffg} we have the following.\vspace*{-2pt}
\begin{proposition}\label{propcloslimcom} Let $X$ be compact and let
$(\rho_1^{(n)},\rho_2^{(n)})$
be a sequence of realizable pairs of symmetric Radon measures which
converges in the vague topology to $(\rho_1,\rho_2)$ and
is such that the
condition of Theorem \ref{propsuffg} holds for $(\rho_1^{(n)},\rho
_2^{(n)})$
for some $R_n \geq0$ with $\liminf_{n \rightarrow\infty} R_n <
\infty$.
Then $(\rho_1,\rho_2)$ is also realizable.\vspace*{-2pt}
\end{proposition}

The proof is similar to the proof of the next theorem, which arises from
Theorem~\ref{chisuff}.\vspace*{-2pt}
\begin{proposition}\label{propcloslimgen} Let $(\rho_1^{(n)},\rho_2^{(n)})$
be a sequence of realizable pairs of symmetric Radon measures which
converges in the vague topology to $(\rho_1,\rho_2)$ and
is such that the
condition of Theorem \ref{chisuff} holds for $(\rho_1^{(n)},\rho_2^{(n)})$
for some fixed $\chi$ and $R_n \geq0$ with $\liminf_{n \rightarrow
\infty} R_n < \infty$.
Then $(\rho_1,\rho_2)$ is also realizable.
\end{proposition}
\begin{pf} We will show that $(\rho_1,\rho_2)$
fulfills the sufficiency condition of Theorem \ref{chisuff}. Without loss
of generality we may replace $(\rho_1^{(n)},\rho_2^{(n)})$ by a~subsequence
such that $R_n$ converges to a finite limit $R$. If
$Q^\chi_{f_0,f_1,f_2,f_3}$ is a~nonnegative $\chi$-restricted polynomial
then the hypotheses imply that
%
%
\begin{equation}
f_0 + \int_X f_1(\x) \rho^{(n)}_1(d\x) + \int_{X \times X}
f_2(\x,\y) \rho^{(n)}_2(d\x,d\y) +f_3 R_n \geq0
\end{equation}
for all $n$. Taking the $n\to\infty$ limit then gives (\ref{chicond}).
\end{pf}

It is easy to see that the conditions of Proposition
\ref{propcloslimgen} may be replaced by the requirement that the pairs
$(\rho_1^{(n)},\rho_2^{(n)})$ can be\vspace*{1pt} realized by
processes~$\mu_n$ in such a way that $\liminf_{n\rightarrow\infty}
\int_{\mathcal{N}(X)}H_3^\chi(\eta)\mu_n(d\eta)<\infty$.

There is an analogous consequence of Theorem \ref{propnc}
whose statement we omit.

\section{Realizability for stationary processes}\label{homogeneous}

In this section we use a variant of the previous results to consider the
question of whether correlation functions having some symmetry can be
realized by a point process having the same symmetry. Throughout this
section we take $G$ to be a topological group acting transitively on
$X$ in
such a way that the action, considered as a map $G \times X \rightarrow X$,
is continuous. The group action can then be extended to an action on
the Radon measures on $X$ and hence, on $\mathcal{N}(X)$ and thus
finally to an
action on point processes; the latter is continuous and linear. We call a
point process \textit{stationary} if it is invariant under this action.
A~stationary point process has stationary correlation functions, that is,
these functions are also invariant under the action of the group. Here we
address the converse question of whether or not stationary correlation
functions can be realized by stationary point processes.

For simplicity we will consider only the possibilities that $G$ be Abelian
or compact, or a semi-direct product of an Abelian and a compact group.

Typical cases are $X=\RR^d$, $\ZZ^d$, etc. As described earlier, there
is then a~natural action of
the translation group on $X$. In this context for a stationary point
process there necessarily exists a real number
$\rho$ such that $\rho_1(d\mathbf{r}) = \rho \,d_I\mathbf{r}$, where
$d_I\mathbf{r}$ denotes the
invariant measure on $X$: Lebesgue measure on~$\RR^d$ and the torus and
counting measure on $\ZZ^d$ and the discrete torus. In general, however,
it may not be true that $\rho_2$ has a density with respect to the Lebesgue
measure on $X^2$; for example, consider on $\mathbb{R}$ the point process
defined by $\mu(d\eta):= \int_0^1 \delta_{\bar{\eta}_y}(d\eta)\,
dy$, where
$\bar{\eta}_y(dr):= \sum_{x \in\ZZ} \delta_{y+x}(dr)$. However,
one can
show that there must exist a Radon measure $g_2$ on $X=\RR^d$ such
that for
any $f_2\in\C_c(X^2)$,
%
%
\begin{equation}\label{red2pt}
\int_X
\int_X f_2(\mathbf{r}_1,\mathbf{r}_2) \rho_2(d\mathbf{r}_1,d\mathbf{r}_2) = \int
_X \int_X f_2(\mathbf{r},\mathbf{r}+
\bar{\mathbf{r}}) \rho^2
g_2(d\bar{\mathbf{r}})\,d_I\mathbf{r}.
\end{equation}
The form in which we have written the right-hand side of
(\ref{red2pt}), isolating a~factor of $\rho^2$ in the
two-point function, is natural in certain applications (see, e.g.,
\cite{CKLS,KLS}).

We first consider the case in which $G$ is Abelian. To be concrete we will
fix a strictly positive bounded continuous function $\chi$ on $X$ and
work in the spirit of Theorem \ref{chisuff}, considering processes with
finite third $\chi$-moments but similar results could be given in the
spirit of Theorem \ref{propnc}. The key idea is to work with processes
satisfying a bound on the third $\chi$-moment which is uniform under the
group action. More precisely, denoting by $g\chi$ the transformed function
$\chi(g\cdot)$, we require a bound for $H_3^{g\chi}$ uniform in $g$.
Proposition \ref{proptrans} establishes the existence of a stationary
process given the existence of one process satisfying such a uniform
bound and Theorem \ref{theosuffgg} gives sufficient conditions, solely in
terms of the given moments, for the realizability by a process satisfying
such a bound.
\begin{proposition}\label{proptrans} Let $G$ be Abelian and let
$(\rho_1,\rho_2)$ be stationary correlation functions realizable
by a process $\mu$ satisfying $\sup_{g \in G}\EE_\mu H_3^{g\chi}
\le R$.
Then $(\rho_1,\rho_2)$ can be realized by a stationary point process.
\end{proposition}
\begin{pf} Let $K_R$ denote the set of all point processes $\mu$ which
realize $(\rho_1,\rho_2)$ and satisfy
$\sup_{g \in G}\EE_\mu H_3^{g\chi} \le R$; $K_R$ is nonempty by
hypothesis. The action of $G$ on point processes leaves $K_R$ invariant.
In Lemma \ref{lemrelcom} we prove that $K_R$ is convex and compact.
Then by
the Markov--Kakutani fixed point theorem (see, e.g., \cite{RSI},
Theorem V.20) there exists a $\mu\in K_R$ which is invariant with
respect to
the action of $G$.
\end{pf}
\begin{lemma}\label{lemrelcom} The set $K_R$ introduced in the proof of
Proposition \ref{proptrans} \textit{is convex and compact in the weak topology}.
\end{lemma}
\begin{pf} The convexity of $K_R$ is obvious. To show that $K_R$ is
compact in the weak topology, we first show that it is tight and hence
precompact. From Lem\-ma~\ref{lemcomp} it follows easily that
$S_N:=\{\eta\in\N(X)\mid\langle\eta,\chi\rangle\le N\}$ is
compact and if
$\mu\in K_R$ then from $\EE_\mu H_3^\chi\le R$ and (\ref{ineq}) it follows
via Markov's inequality that for $\epsilon>0$ there is a choice of $N$,
depending only on $\epsilon$ and $R$, such that
$\mu(S_N) >1-\epsilon$, verifying tightness.

It remains to prove that $K_R$ is closed. Let $\mu_n$ be a sequence in
$K_R$ which converges weakly to a point process $\mu$. Approximating
$H_3^{g\chi}$ by an increasing sequence of bounded continuous functions
and using the convergence of the sequence $\mu_n$ on such functions
and the
monotone convergence theorem for $\mu$, we find that
$\int_{\mathcal{N}(X)}H_3^{g\chi}(\eta) \mu(d\eta) \leq R$. It
remains to
show that $\mu_n$ converges also on every quadratic polynomial
$P=P_{f_0,f_1,f_2}$ with $f_1\in\C_c(X) $ and $f_2\in\C_c(X^2)$, which
guarantees that $\mu$ has the correct first and second correlation
functions. But by~(\ref{ineq}), $|P(\eta)|\le A+B\E{\chi,\eta}^2$
for some
$A,B\ge0$ and so for $M\ge2A$, $|P(\eta)|\ge M$ implies
$|P(\eta)|\le2B\E{\chi,\eta}^2$ and so for any $\nu\in K_R$,
%
%
\begin{eqnarray}\label{cutoff}
\int_{P\ge M} |P(\eta)| \nu(d\eta)
&\le& 2B\int_{\E{\chi,\eta}^2\ge M/2B}\E{\chi,\eta}^2 \nu
(d\eta)
\nonumber\\
&\leq& \frac{(2B)^{3/2}}{M^{1/2}}\int_X\E{\chi,\eta}^3 \nu
(d\eta)\\
&\leq&\frac{(2B)^{3/2}}{M^{1/2}}2(b_3+R),\nonumber
\end{eqnarray}
where we have used (\ref{ineq}) again. But if
$P^{(M)}(\eta):=\operatorname{sign}[P(\eta)]\min\{|P(\eta)|,M\}$
then for
any fixed $M$,
\[
\int_XP^{(M)}(\eta) \mu_n(d\eta) \longrightarrow
\int_XP^{(M)}(\eta) \mu(d\eta) \qquad\mbox{as $n\to\infty$}
\]
and with (\ref{cutoff}) the proof is complete.
\end{pf}

Our sufficient condition for the existence of a process, analogous to
Theorem~\ref{chisuff}, involves polynomials of the form
%
%
\begin{eqnarray}\label{eqadmg}
&&Q^\chi_{f_0,f_1,f_2,(f_{3,1},g_1),\ldots,(f_{3,n},g_n)}(\eta)\nonumber\\[-8pt]\\[-8pt]
&&\qquad
= f_0 + \langle f_1, \eta\rangle+ \langle f_2 ,
\eta^{\odot2} \rangle+ \sum_{i=1}^n f_{3,i}
H^{g_i\chi}_3(\eta),\nonumber
\end{eqnarray}
where $\chi$ is as above, $f_0$ and
$f_{3,1},\ldots,f_{3,n}$ are real numbers, $f_1$ and $f_2$ are continuous
symmetric functions with compact support on $X$ and $X\times X$,
respectively, and $g_1,\ldots,g_n\in G$. The term
$\sum_{i=1}^n f_{3,i} H^{g_i\chi}_3$ in (\ref{eqadmg}) controls moments
involving $H^{g\chi}_3$ and also makes the set of all the $Q^\chi$
into a
vector space.
\begin{theorem}\label{theosuffgg} Let $G$ be Abelian and let $\rho_1$
and $\rho_2$ be symmetric $G$-sta\-tionary Radon measures on $X$ and
$X\times X$, respectively. Then $\rho_1$
and $\rho_2$ are realizable by a stationary point process
$\mu$ with $\sup_{g \in G} \int H^{g \chi}_3(\eta) \mu(d\eta) <
\infty$
if and only if there is a constant $R>0$
such that if
$Q^\chi_{f_0,f_1,f_2,(f_{3,1},g_1),\ldots,(f_{3,n},g_n)}$ is nonnegative
on $\N(X)$ then
%
\begin{equation}\quad
f_0 + \int_X f_1(\x) \rho_1(d\x) + \int_{X \times X} f_2(\x
,\y)
\rho_2(d\x,d\y) +\sum_{i=1}^n f_{3,i} R \geq0.
\end{equation}
\end{theorem}
\begin{pf} The proof is completely parallel to the proofs of
Theorems \ref{propsuffg} and~\ref{chisuff} and we mention only a few
details. Let $V$ be the vector space of all functions which have the form
$F + \sum_{i=1}^n \alpha_i H_3^{g_i \chi}$, where $F$ is a continuous
function on $\mathcal{N}(X)$ satisfying $|F|\le C(1+H_3^\chi)$ for some
constant $C>0$, $\alpha_1,\ldots,\alpha_n \in\RR$ and
$g_1, \ldots,g_n \in G$. Let $V_0$ be the subspace of $V$ consisting of
all polynomials $Q^\chi_{f_0,f_1,f_2,(f_{3,1},g_1),\ldots,(f_{3,n},g_n)}$.
For any $R\ge0$ we define a linear functional $L_R$ on $V_0$ by
\begin{eqnarray*}
&&L_R\bigl(Q_{f_0,f_1,f_2,(f_{3,1},g_1), \ldots, (f_{3,n},g_n)}\bigr)\\
&&\qquad
:= f_0 + \int_X f_1(\mathbf{r}) \rho_1(d\mathbf{r}) + \int_{X \times X}
f_2(\mathbf{r}_1,\mathbf{r}_2) \rho_2(d\mathbf{r}_1,d\mathbf{r}_2)+ \sum
_{i=1}^n f_{3,i} R
\end{eqnarray*}
and show that $\rho_1,\rho_2$ is realizable by a process $\mu$ with
$\sup_{g \in G} \int H^{g \chi}_3(\eta) \mu(d\eta) < \infty$ if
and only if
$L_R$ is nonnegative for some $R>0$. The condition is clearly necessary.
Conversely, if $R$ is such that $L_R$ is nonnegative on $V_0$, we may extend
$L_R$ to $V$ using Theorem \ref{RKET}. To show that this extended linear
form is given by a measure we let $W$ be the\vspace*{1pt} subspace of $V$ consisting of
all continuous functions $F\in V$ such that $|F|\le C(1+H_2^\chi)$ for some
$C>0$ and apply Theorem \ref{theodaniel} to $L_R$ on $W$. The verification
that $L_R$ satisfies (D) on $W$ is the same as the corresponding
verification for Theorem \ref{chisuff} and we conclude that there
exists a
probability measure $\mu$ on $\mathcal{N}(X)$ such that for $F\in W$,
\[
L_R(F) = \int_{\mathcal{N}(X)} F(\eta) \mu(d\eta).
\]
As $W$ includes all $Q_{f_0,f_1,f_2}$ the measure $\mu$ realizes
$(\rho_1,\rho_2)$. Finally, for $n\in\NN$ and $g \in G$ the lower
semi-continuous function $H^{g\chi}_3$ can be approximated from below
by an
increasing sequence of continuous bounded functions~$H^{g\chi}_{3,k}$. By
the positivity of $L_R$ on $V$,
$\int
H^{g\chi}_{3,k} \,d\mu=L_R(H^{g\chi}_{3,k})\le L_R(H^{g\chi}_3)=
L_R(Q_{0,0,0,(1,g)})=
R$ and so the monotone convergence theorem implies that $\int
H^{g\chi}_3 \,d\mu\le R$. The result follows from Proposition
\ref{proptrans}.
\end{pf}

Next we consider the case of compact groups.
\begin{proposition}\label{prosuffgna} Let $G$ be a compact and let
$\rho_1$ and $\rho_2$ be symmetric $G$-stationary Radon measures on
$X$ and
$X\times X$. Then $\rho_1$ and $\rho_2$ are realizable by a stationary
point process $\mu$ if and only if they are realizable.
\end{proposition}
\begin{pf} Let $\mu$ be a realizing point process for $\rho_1$ and
$\rho_2$. Denote by $\nu$ the Haar measure on $G$ and by $g\mu$ the point
process transformed via the action of $g$. Then define
$\tilde{\mu} := \int_G (g\mu) \nu(dg)$ in the sense that
\[
\int F(g\eta)\tilde{\mu}(d\eta) := \int_G F(g \eta)\mu(d\eta)
\nu(dg) \qquad\mbox{for all } F \in L^{1}(\N(X),\mu),
\]
$\tilde{\mu}$ is a stationary realizing point process.
\end{pf}

Finally, we may easily combine the previous two cases and, in particular,
cover the important special case of the Euclidean group acting on $\RR^n$.
\begin{proposition}\label{corsemi} Let $G$ be the semi-direct product
$N\rtimes H$ of
an Abelian group $N$ and a compact topological group $H$ and let $\rho_1$
and $\rho_2$ be symmetric $G$-stationary Radon measures on $X$ and
$X\times X$. Then $\rho_1$ and $\rho_2$ are realizable by a stationary
point process $\mu$ with
$\sup_{g \in G} \int H^{g \chi}_3(\eta) \mu(d\eta) < \infty$ if
and only if
there exists a constant $R>0$ such that if
$Q^\chi_{f_0,f_1,f_2,(f_{3,1},g_1),\ldots,(f_{3,n},g_n)}$, $g_i \in
N$, is
nonnegative on $\N(X)$ then
%
%
\begin{equation}\quad
f_0 + \int_X f_1(\x) \rho_1(d\x) + \int_{X \times X} f_2(\x
,\y)
\rho_2(d\x,d\y) +\sum_{i=1}^n f_{3,i} R \geq0.
\end{equation}
\end{proposition}
\begin{pf} Applying Theorem \ref{theosuffgg} to the action of $N$ we
obtain an $N$-sta\-tionary point process. Using the construction in
Proposition \ref{prosuffgna} we arrive at point process also stationary
under the action of $H$ and hence, stationary for the action of $G$.
The particular structure of the multiplication in the semi-direct product
does not play any role.
\end{pf}

As a closing remark, we note that in this section we have concentrated on
extensions of the results of Section \ref{subsubnoncompact} to
stationary processes but that corresponding extensions for the results in
Sections \ref{susechard} and \ref{compact} can be obtained
similarly and
in fact more easily. The next proposition gives extensions of
Theorems \ref{thsuflat} and \ref{propsuffg}.
\begin{proposition}
Let $G$ be as in Proposition
\ref{corsemi} and let $\rho_1$
and $\rho_2$ be symmetric $G$-stationary Radon measures on $X$ and
$X\times X$. Then:

\textup{(a)} If $\rho_2$ forces a hard core exclusion for a metric $d$
and the
action of~$G$ leaves $d$ invariant, then $\rho_1$ and $\rho_2$ are
realizable by a stationary point process~$\mu$ if and only if they are
realizable by a point process.

\textup{(b)} If $X$ is compact, then $\rho_1$ and $\rho_2$ are
realizable by a
stationary point process $\mu$ with finite third moment if and only if they
are realizable by a point process with finite third moment.
\end{proposition}
\begin{pf} In each case one first verifies the result for $G$ Abelian
and then extends to the semi-direct product case as in the proof of
Proposition \ref{corsemi}. When $G$ is Abelian the proof of (a) follows
the proof of Proposition \ref{proptrans} but now instead of $K_R$ we
consider the set $K$ of all measures realizing $(\rho_1,\rho_2)$. $K$~is
obviously convex; to show that $K$ is compact we note that since $\N_D(X)$
is compact so is the set of all probability measures on $\N_D(X)$~%
\cite{P67}, of which~$K$ is a subset. Moreover, $K$ is closed since,
because quadratic polynomials on~$\N_D(X)$ are bounded and continuous, weak
limit points of $K$ give the same expectation values of quadratic
polynomials as do points in~$K$ and thus also realize
$(\rho_1,\rho_2)$. Part (b) follows from Proposition \ref{proptrans}
itself by taking $\chi=1$ there and using the fact that then $g\chi
=\chi$
for all $g\in G$.
\end{pf}


\section{Classes of necessary conditions}\label{secclass}

If $X$ is finite then, as indicated in Remark \ref{hcremark}(c), the
necessary and sufficient conditions of Theorem \ref{thsuflat} give
rise to
a finite linear programming problem. In this section we allow $X$ to be
infinite and consider the problem of isolating useful necessary conditions
from among the uncountably infinite class of Theorem \ref{lemnes}. In
the latter case, as in the former, the conditions arising from distinct
functions may be related; in particular, some of them may imply others. For
practical purposes it would be desirable to identify a class of functions,
as small as possible, such that the conditions arising from this class
imply all the conditions but for this presumably very hard problem we have
no solution at the moment. In this section we will, however, for a certain
uncountable subclass of the full class of conditions of
Theorem \ref{lemnes}, identify a handful of conditions which imply
those of
the whole subclass so that one may check all conditions arising from the
subclass by checking the few selected conditions.

Suppose that we are given a pair $(\rho_1,\rho_2)$ of correlation functions
and wish to use Theorem \ref{lemnes} to show that this pair is not
realizable on some $\Ns$. If~$\rho_2$ forces a hard core exclusion with
diameter $D$ or if we impose a bound on the number of particles as in
Corollary \ref{CoNrep}, then we would take $\Ns$ to be $\N_D(X)$,
$\N^{\leq N}(X)$ or $\N^N(X)$, but otherwise we have a priori no
better choice than to take $\Ns=\N(X)$. The general strategy that we
suggest, and will illustrate by an example, is to introduce a family of
polynomials on~$\N(X)$ depending on some finite set of parameters and then
to determine a~finite subset of this family such that satisfying the
necessary conditions for polynomials in the subset guarantees satisfaction
for all polynomials in the original family.

We work out this strategy in a particular case obtaining in the process
several standard necessary conditions which have appeared in the
literature (see \cite{KLS} for a detailed exposition and references).
We choose a fixed nonzero $f\in\M_c(X)$ and consider the
family of all polynomials of the form
%
%
\begin{equation}
P^{(a,b,c)}(\eta):=a\langle f,\eta\rangle^2+b\langle f,\eta\rangle+c.
\end{equation}
Note that in the notation of (\ref{quadpoly}),
$P^{(a,b,c)}=P_{f_0,f_1,f_2}$, with
\[
f_2(\vect{r}_1,\vect{r}_2)= a f(\vect{r}_1) f(\vect{r}_2),\qquad
f_1(\vect{r}) = b f(\vect{r})+ a f^2(\vect{r}),\qquad
f_0 = c.
\]
Let
$F:=\{\langle
f,\eta\rangle\mid\eta\in\mathcal{N}_{\mathrm{supp}}\}\subset\RR$
be the range of $\langle f,\cdot\rangle$ and let $\Gamma$ be the convex
cone of all $(a,b,c) \in\RR^3$ such that $p(x):=ax^2+bx+c \geq0$ for all
$x \in F$. The necessary condition then is that the linear function
$L=L_{\rho_1,\rho_2}$ defined by
%
%
\begin{eqnarray}
\label{neccondL}
L(a,b,c)&=& a \int_{X^2} f(\mathbf{r}_1) f(\mathbf{r}_2)
\rho_2(d\mathbf{r}_1,d\mathbf{r}_2)\nonumber\\[-10pt]\\[-10pt]
&&{} + \int_X \bigl(b f(\mathbf{r}) + a f^2(\mathbf{r})
\bigr) \rho_1(d\mathbf{r}) + c,\nonumber
\end{eqnarray}
should be nonnegative on $\Gamma$.\vadjust{\eject}

Before continuing we give a (nonexhaustive) discussion of possible
structure of $F$, excluding the uninteresting case $f=0$, in order to give
some feeling for how this structure can affect the necessary conditions.
If $\Ns=\N_D(X)$ then $F$ is bounded above and below, for example, by
$\pm M_D\sup|f|$, where~$M_D$ is the maximum
number of disjoint balls of diameter $D$ which can be placed so that their
centers lie in the support of $f$. Similar bounds hold if~$\Ns$ is
$\N^{\leq N}(X)$ or $\N^N(X)$. Otherwise $F$ is unbounded and is bounded
below (by~0) if and only if $f\ge0$ and above (again by 0) if and only if
$f\le0$. If~$f$ takes only a finite number of values then $F$ will consist
of certain linear combinations, with integer coefficients, of these values
and $F$ may then be discrete or may be dense in $\RR$; in the simplest
case, when $f=\Ii_\Lambda$ for some $\Lambda\subset X$, $F$ is just
a set
of nonnegative integers. If $\Ns=\N(X)$, $f$ is nonnegative and the
range of
$f$ contains some interval $(0,\delta)$, then $F=\RR_+$, or if the
range of
$f$ contains some interval $(-\delta,\delta)$, then $F=\RR$.

We make two more preliminary remarks. First, if a realizing measure~%
$\mu$ exists then $E(f) := \EE_\mu\E{f,\cdot}$ and
$V(f) := \Var_\mu(\E{f,\cdot})$ may be calculated from~$\rho_1$ and
$\rho_2$ as
\begin{eqnarray*}
E(f) &=& \int_X f(\mathbf{r}) \rho_1(d\mathbf{r}),\\
V(f) &=&\int_{X^2} f(\vect{r}_1) f(\vect{r}_2)
\rho_2(d\vect{r}_1,d\vect{r}_2) \\
&&{}
+ \int_X f(\vect{r})^2 \rho_1(d\mathbf{r})-
\biggl( \int_X f(\vect{r}) \rho_1(d\vect{r}) \biggr)^2,
\end{eqnarray*}
so that
%
%
\begin{equation}\label{simple}
L(a,b,c)=a V(f) + p[E(f)].
\end{equation}
Second, due to the homogeneity in $(a,b,c)$ of the problem it suffices to
consider conditions arising from polynomials with either $a=0$ or
$a=\pm1$.

\textit{Case} 1. $a=1$. In this case, (\ref{simple}) implies that the
constraint on $\rho_1,\rho_2$ will be of the form $V(f)\ge-p[E(f)]$; by
taking $p(x)=[x-E(f)]^2$ we recover the obvious requirement that
$V(f)\ge0$. The condition that $V(f)\ge0$ for all $f \in\C_c(X)$
is equivalent to the so-called variance condition; cf., for example,~%
\cite{KLS}.
If $E(f)\in F$ then $p[E(f)]\ge0$ whenever $p\in\Gamma$ so
that (\ref{simple}) implies that for no choice of $b$ and $c$ can
$L(1,b,c)\ge0$ impose further restrictions on $\rho_1,\rho_2$. Otherwise,
$E(f)\in(x_-,x_+)$ for some maximal open interval $(x_-,x_+)$ disjoint from
$F$; then the choice $p_0(x)=(x-x_-)(x-x_+)$ implies the constraint
%
%
\begin{eqnarray}\label{age0}
V(f)&\ge&\bigl(x_+-E(f)\bigr)\bigl(E(f)-x_-\bigr) \hspace*{120pt}\nonumber\\[-8pt]\\[-8pt]
&&\eqntext{\mbox{for }x_- := \sup\{x \in F\mid x \leq E(f)\},\qquad
x_+ := \inf\{x \in F\mid x \geq E(f)\}.}
\end{eqnarray}
An easy computation shows that any monic quadratic polynomial $p$ with
$p(x_-)$, $p(x_+)\ge0$ satisfies $p[E(f)]\ge p_0[E(F)]$, so that (\ref{age0})
includes all restrictions arising in Case 1 [note that as written the
constraint (\ref{age0}) includes the case $E(f)\in F$]. If $f=\Ii
_\Lambda$
for $\Lambda\subset X$ then $F=\NN_0$ and (\ref{age0}) was found by Yamada
\cite{Yamada}. Whether for other choices of $f$ one obtains additional
restrictions is unknown.

In the case $x_-,x_+ \in F$, $x_- < x_+$
the choice of $x_-,x_+$ in (\ref{age0})
corresponds to an extremal ray in the cone $\Gamma$. The cone can be
defined as intersection $\bigcap_{y \in F} H_y$ with
$H_y :=\{ (a,b,c) \in\RR^3\mid ay^2+by+c\geq0 \}$. Hence, to each pair
$x_1\leq x_2 \in F$ there corresponds a ray
$\{(a,b,c)
\in\RR^3\mid ax_1^2+bx_1+c = 0$ and $ax_2^2+bx_2+c = 0
\}$.
This ray will be in the cone and hence, an extremal ray only if
$(x_1,x_2) \cap F = \varnothing$. Hence, the choice of $x_-,x_+$ in
(\ref{age0}) corresponds to a particular extremal ray of $\Gamma$.

\textit{Case} 2. $a=-1$. In this case, $p(x)$ can be nonnegative on $F$ only
if $F$ is bounded and reasoning as in the previous case shows that the
constraint obtained from $p(x)=(\sup F-x)(x-\inf F)$,
%
%
\begin{equation}\label{ale0}
V(f)\le[\sup F-E(f)][E(f)-\inf F]
\end{equation}
implies all others.

\textit{Case} 3. $a=0$. We assume $b\ne0$ since a constant polynomial
conveys no restriction; then we may take $b=\pm1$ and thus consider
$p(x)=\pm(x-x_0)$. Such a linear function can be nonnegative on $F$ only
if either (i) $F$ is bounded below, in which case the constraint from
$p(x)=x-\inf F$ implies all others, or (ii) $F$ is bounded above, in which
case a similar conclusion holds for $p(x)=\sup F-x$. If $f$ is nonnegative
then $\inf F=0$ and the condition in (i), $E(f)\ge0$, just asserts the
positivity of the measure $\rho_1$. Case (ii), namely, $E(f)\leq\sup F$,
can occur if $\rho_2$
enforces a hard core exclusion or if we impose an a priori bound on
the number of particles as in Corollary \ref{CoNrep}. A simple
interpretation can be given when $X$ is compact and $f=\Ii_X$; then
$\sup F$ is the maximum number $M$ of points which can be contained in $X$
under the hard core or {a priori} condition and the condition imposed
on $\rho_1$ by the constraint $E(\Ii_X)=\int_X\rho_1(dx)\le M$ is
that the
expected number of points be less than this maximum. If we further assume
that $X$ is a torus with Lebesgue measure $\nu$ and that
$\rho_1(dx)\equiv\rho\nu(dx)$ is invariant under translations
then this
condition is $\rho\le M/\nu(X)$. We can then see that no constraint
arising from another choice of $f$ in Case 3 gives
further restrictions on $\rho_1,\rho_2$; indeed, since picking an
$\eta\in\Ns(X)$ with $\eta(X)=M$ (which one does not matter) we
find from
$\rho\le M/\nu(X)$ that for any $f$
%
%
\begin{equation}
E(f)=\rho\int_Xf(x) \nu(dx)
=\rho\int_X\frac1{M}\sum_{y\in\eta}f(y+x) \,d\nu(x)
\le\sup F,
\end{equation}
because $\sum_{y\in\eta}f(y+x)
\le\sup F$.

\section*{Acknowledgments}

We are grateful to J. K. Percus for fruitful discussions and providing
us with references and unpublished material. We would like to thank
V.~Bach, G. Moreano and M. Esguerra for pointing out references in
quantum chemistry and two anonymous referees for helpful comments. We
would like to thank IHES, Paris for their hospitality.
T. Kuna would like to thank A.v.H. for support and
Yu. G. Kondratiev for discussion and references.


%
\printaddresses

\end{document}